\documentclass{article}
\usepackage{a4}
\title{Sheaves for predicative toposes}
\author{Benno van den Berg}
\date{22 July, 2005}
\usepackage{amssymb}
\usepackage{amsmath}

\newtheorem{proposition}{Proposition}[section]
\newtheorem{lemma}[proposition]{Lemma}
\newtheorem{corollary}[proposition]{Corollary}
\newtheorem{definition}[proposition]{Definition}
\newtheorem{theorem}[proposition]{Theorem}

\newtheorem{rem}[proposition]{Remark}
\newcounter{opgaveteller}
\newcounter{lijst-teller}

\newcounter{boon}
\newcounter{boon2}
\newenvironment{proof}{\noindent {\bf Proof:} \nopagebreak }{\nopagebreak\hfill$\Box$}

\newenvironment{remark}{\begin{rem} \rm}{\end{rem}}

\input xypic
\xyoption{all}
\UseComputerModernTips
\CompileMatrices
\definemorphism{dat}\dashed\tip\notip
\definemorphism{dub}\Solid\Tip\notip

   {\begin{list}{\roman{lijst-teller}\/)\hfil}%
              {\labelwidth 2em%
               \leftmargin\labelwidth\advance\leftmargin by\labelsep%
               \usecounter{lijst-teller}}}%
   {\end{list}}

   {\end{list}}

    {\end{list}}

\newdir{|>}{!/2pt/@{|}*@^{>}*@_{>}}
\newdir{ >}{{}*!/-5pt/@{>}}
\newdir{ )}{{}*!/-5pt/@{)}}

\begin{document}

\maketitle

\begin{quote}
{\bf Abstract:}  In this paper, we identify some categorical structures in which one can model predicative formal systems: in other words, predicative analogues of the notion of a topos, with the aim of using sheaf models to interprete predicative formal systems. Among our technical results, we prove that all the notions of a ``predicative topos'' that we consider, are stable under presheaves, while most are stable under sheaves.
\end{quote}

\section{Introduction}

The importance for topos theory to logic is due to the fact that toposes provide categorical models of impredicative constructive theories, like higher-order type theoretic logics or the set theory IZF. The theory of toposes yields a large stock of examples of such models in the form of toposes of sheaves, thereby also incorporating the set-theoretic method of forcing. Toposes of sheaves for a site are an especially fruitful source of examples, because the construction can be iterated: i.e., the notion of a sheaf for a site can be formulated internally in a topos and these sheaves for the internal site again form a topos. Using these topos-theoretic models, one can obtain consistency and independence results and derived rules (good sources for logical applications of topos theory are \cite{BlassScedrov89}, \cite{TroelstraVanDalen88}, Chapter 15, \cite{MacLMoe:Sheaves}, Chapter 6)

But important instances of constructive formal theories are also predicative. One could think of Myhill's set theory CST, extended by Peter Aczel to CZF and Martin-L\"of type theory, where CST and CZF have the predicative features of restricting the use of the separation and powerset axiom.

For the study of metamathematical properties of predicative formal systems, topos-theoretic methods, like sheaves, have been much less developed, let alone used. As it happens, the notion of a topos is not the most suitable for that purpose, as it is an inherently impredicative structure by containing power objects and a subobject classifier. So when one tries to formulate an appropriate notion of sheaf for a predicative formal system, one first has to modify the notion of a topos. Therefore in their work \cite{MoerPalm:Wellftr} and \cite{MoerPalm:CST}, Moerdijk and Palmgren start to investigate predicative analogues of the notion of a topos and to redevelop for these ``predicative toposes'' those parts of the theory of toposes that are important for logical applications, like the theory of sheaves.

The theory of predicative toposes is connected to type theory and set theory in the following ways. The connection of predicative toposes with type theory relies on the work of Seely \cite{Seely84}, who established a correspondence between locally cartesian closed categories and type theories (for coherence problems related to substitution, people have found different solutions, for example \cite{Hofmann94}). In \cite{MoerPalm:CST}, Moerdijk and Palmgren establish the connection with set theory by building a model of Aczel's CZF in a sufficiently strong version of a predicative topos. Because of these relations, the theory of sheaves in the context of predicative toposes has implications for type theory and constructive set theory.

The work of Moerdijk and Palmgren provide a first step in formulating the right notion of a predicative topos. They suggest the notions of a $\Pi W$-pretopos and a stratified pseudotopos and in this paper, we propose several other related axiomatisations. The candidates are introduced in the second section of this paper. The main parts of the paper are sections 4 and 5 in which we show that all candidates for the title of predicative topos that we consider are stable under presheaves, while most are also stable under taking sheaves for an internal site. We hope to take up the logical implications of these results in later work.

More concretely, we show that:
\begin{enumerate}
\item When $\bf E$ is a $\Pi W$-pretopos and $\bf C$ is an internal category, then PSh$_{\bf E}({\bf C})$ is also a $\Pi W$-pretopos. In \cite{MoerPalm:Wellftr}, Moerdijk and Palmgren assume an additional principle of transfinite induction for W-types, which is shown to be superfluous.
\item When $\bf E$ is a $\Pi W$-pretopos with a class of small maps satisfying the collection axiom and $\bf C$ is an internal site with small covers, then Sh$_{\bf E}({\bf C})$ is a $\Pi W$-pretopos. In \cite{MoerPalm:CST}, the authors make the assumption that the small maps satisfy an additional axiom, the axiom of multiple choice (AMC). We show that this assumption is not necessary.
\item We propose an alternative notion of predicative topos, weaker than Moer\-dijk and Palmgren's notion of a stratified pseudotopos (with AMC): a $\Pi W$-pretopos satisfying the universe operator axiom UO. We show that for a $\Pi W$-pretopos $\bf E$ satisfying UO and containing an internal site $\bf C$, the category Sh$_{\bf E}({\bf C})$ is a $\Pi W$-pretopos again satisfying UO. A similar stability result holds for the stratified pseudotoposes, as shown in \cite{MoerPalm:CST}, but the proof for our alternative is much simpler. 
\end{enumerate}

Besides Moerdijk and Palmgren, other people have worked on sheaves for predicative formal systems. An older source is the work by Grayson \cite{Grayson83}, while some recent work has been done by Gambino \cite{Gambino02}. Their work is done in the context of set theory, while some (unpublished) work in connection with type theory has been done by Martin-L\"of. Awodey and Warren have a categorical approach in \cite{AwodeyWarren05}, but are not concerned with the same issues as we are here.

I would like to thank several people for discussing with me the research reported here. First of all, my supervisor Ieke Moerdijk, but also Martin Hyland, Nicola Gambino, Thierry Coquand, Per Martin-L\"of and Erik Palmgren, who also invited me to speak on these matters in the Stockholm-Uppsala Logic Seminar. I also spoke on this topic during the Summer School on Topos Theory in Haute-Bodeux, and I would also like to thank the organizers of the summer school for giving me this opportunity.

\section{A predicative notion of topos}

\subsection{Categorical preliminaries}

We use this section to introduce some terminology and recall some useful facts from category theory.

\begin{definition}
A subobject $X$ of $A$ in a category $\bf E$ is called the \emph{image} of $f: B \to A$, if $X$ is the least subobject through which $f$ factors. A map $f: B \to A$ in $\bf E$ is called a \emph{cover}, when its image is the maximal subobject $A$ of $A$.
\end{definition}
\begin{definition}
A square
\begin{displaymath}
\xymatrix{ D \ar[r] \ar[d] & B \ar[d] \\
           C \ar[r]        & A }
\end{displaymath}
in a category $\bf E$ with pullbacks is called a \emph{quasi-pullback}, if the induced map \[ \xymatrix{D  \ar[r] & B \times_{A} C} \] is a cover.
\end{definition}
\begin{definition}
A category $\bf E$ is called \emph{regular}, when it has finite limits and all maps factor, in a stable fashion, as a cover followed by a mono. A regular category $\bf E$ is called \emph{exact}, when every equivalence relation
\begin{displaymath}
\xymatrix{R \ar@<.7ex>[r] \ar@<-.7ex>[r] & X }
\end{displaymath}
has a quotient $X/R$, i.e., an object that fits into a diagram
\begin{displaymath}
\xymatrix{R \ar@<.7ex>[r] \ar@<-.7ex>[r] & X \ar[r] & X/R }
\end{displaymath}
that is both a pullback and a coequalizer. In an exact category, such quotients are moreover required to be stable under pullback.
\end{definition}
\begin{definition} A category $\bf E$ is called a \emph{pretopos}, when it is exact and has finite disjoint and stable sums.
\end{definition}

In a pretopos, epis and covers coincide.

In general, when $\bf E$ is a category with finite limits in which for any map $f: B \to A$ the induced pullback functor $f^{*}: {\bf E}/A \to {\bf E}/B$ has a right adjoints $\Pi_f$, the category $\bf E$ is called \emph{locally cartesian closed}. A locally cartesian closed pretopos will also be called a \emph{$\Pi$-pretopos}. A $\Pi$-pretopos has enough structure to interpret first-order intuitionistic predicate logic.

In a $\Pi$-pretopos, any map $f: B \to A$ gives rise to its associated polynomial functor:
\[ P_f(X) = \Sigma_{a \in A} X^{B_a} \]
If an initial algebra for this endofunctor exists, it is called the W-type for $f$. We think of $f$ as specifying a signature: a term constructor for every element $a \in A$ of arity $B_a$. The W-type $W_f$ is, whenever it exists, the object of all terms over the signature specified by $f$. Since $W_f$ is a $P_f$-algebra, it naturally comes equipped with a map $P_f(W_f) \to W_f$, usually denoted by sup. Intuitively, for any pair $(a \in A, t: B_a \to W_f)$, sup$_a(t)$ is the term constructed by taking $a$ and substituting $t(b)$ in the $b$-th component of the term $a$.

In the sequel, we need the following characterization theorem (see \cite{BvdB:IndT}):
\begin{theorem}
In a $\Pi$-pretopos $\bf E$ with a natural number object, a $P_f$-algebra $(V, m: P_f(V) \to V)$ is the W-type for a morphism $f$ iff it has no proper $P_f$-subalgebras and $m$ is iso.
\end{theorem}

We will also need the notion of a subterm and the fact that this notion can be formalized in the internal logic of $\Pi$-pretoposes with a natural number object. Briefly, when the W-type $W_f$ of a map $f$ exists, one calls a sequence of the form
\[ \langle w_0, b_0, w_1, b_1, \ldots, w_n \rangle \]
a \emph{path from $w_0$ to $w_n$}, if $w_i \in W_f$ and $b_i \in B$ are such that they satisfy the following compatibility condition: if for an $i<n$, $w_i$ is of the form sup$_{a_i} t_i$, then $f(b_i) = a_i$ and $w_{i+1}=t_i b_i$. Then $v$ is a \emph{subterm} of $w$, when there is a path from $w$ to $v$. (For more details on the formalization of paths and subterms in the internal logic of $\Pi$-pretoposes with natural number object, see \cite{BvdB:IndT}.)

A $\Pi$-pretopos $\bf E$ in which all W-types exist is a \emph{$\Pi W$-pretopos}. A $\Pi W$-pretopos always has a natural number object, because it is the W-type associated to one of the sum inclusions $1 \to 1+1$.

\subsection{Small maps}

Among the various notions of predicative topos that we will discuss in this paper, the concept of a $\Pi W$-pretopos is the most basic. The main problem with this concept is that we will not able to show that $\Pi W$-pretoposes are stable under taking sheaves for an internal site. A natural solution is strenghtening the notion of a predicative topos by formulating a categorical analogue of a type theoretic universe. To this end, we introduce ideas from algebraic set theory. The basic context is that of a category $\bf E$ equipped with a class of maps $S$. The maps in $S$ are referred to as small and the intuition is that their fibres are small in some sense. One should think of: a set as opposed to a proper class, finite as opposed to infinite, countable as opposed to uncountable, a small type as opposed to a type outside the universe of small types, etcetera. The idea is not to fix a set of axioms for this class $S$ once and for all, but algebraic set theory is supposed to be a flexible framework for the categorical study of notions from set theory. Still, we will have to make a choice in order to get started. We will follow \cite{MoerPalm:CST}\footnote{For different axiom systems, see \cite{JoMoe:AST}, \cite{Awodeyenco:AST} and other references at the ``Algebraic Set Theory'' website: {\tt http://www.phil.cmu.edu/projects/ast/}. The main reason for axiomatizing the notion of a class of small maps as we do, is to include the category of setoids as a natural example. Still, the results presented here should be largely independent of such choices.}.

So let $\bf E$ be $\Pi W$-pretopos and let $S$ be a class of maps.
\begin{definition} $S$ is called \emph{stable} if it satisfies the following axioms:
\begin{itemize}
\item[S1] {\rm (Pullback stability)} In a pullback square
\begin{equation} \label{pullb}
\begin{gathered}
\xymatrix{ D \ar[r] \ar[d]_g & C \ar[d]^f \\
           B \ar[r]_p        & A }
\end{gathered}
\end{equation}
$g$ belongs to $S$, whenever $f$ does.
\item[S2] {\rm (Descent)} If in a pullback diagram as in (\ref{pullb}), $p$ is epi, then $f$ belongs to $S$, whenever $g$ does.
\item[S3] {\rm (Sum)} If two maps $f: B \to A$ and $f': B' \to A'$ belong to $S$, then so does $f+f': A+A' \to B+B'$.
\end{itemize}
\end{definition}
These axioms express that maps belong to $S$ in virtue of the properties of their fibres.

\begin{definition} A stable class $S$ is called a \emph{locally full subcategory}, if it also satisfies the following axiom:
\begin{itemize}
\item[S4] In a commuting triangle
\begin{displaymath}
\xymatrix{ C \ar[rr]^g \ar[dr]_h & & B \ar[ld]^f \\
           & A & }
\end{displaymath}
where $f$ belongs to $S$, $g$ belongs to $S$ if and only if $h$ does.
\end{itemize}
\end{definition}

\begin{remark} If S1 holds and all identities belong to $S$, S4 is equivalent to the conjunction of the following two statements:
\begin{itemize}
\item[S4a] Maps in $S$ are closed under composition.
\item[S4b] If $f: X \to Y$ belongs to $S$, the diagonal $X \to X \times_Y X$ in ${\bf E}/Y$ also belongs to $S$.
\end{itemize}
When thinking in terms of type constructors, this means that S4 expresses that smallness is closed under dependent sums and (extensional) equality types. We will actually require the class of small maps to be closed under all type constructors, hence the next definition.
\end{remark}

For any object $X$ in $\bf E$, we write $S_X$ for the full subcategory of ${\bf E}/X$ whose objects belong to $S$. An object $X$ is called \emph{small}, when the unique map $X \to 1$ is small.

\begin{definition}
A locally full subcategory $S$ in a $\Pi W$-pretopos ${\bf E}$ is called a \emph{class of small maps}, if for any object $X$ of ${\bf E}$, $S_X$ is a $\Pi W$-pretopos, and the inclusion functor $S_X \to {\bf E}/X$ preserves the structure of a $\Pi W$-pretopos.
\end{definition}

\begin{definition} A stable class (locally full subcategory, class of small maps) $S$ is called \emph{representable}, if there is a map $\pi: E \to U$ in $S$ such that any map $f: B \to A$ in $S$ fits into a double pullback diagram of the form
\begin{displaymath}
\xymatrix{B \ar[d]_f & B' \ar[r] \ar[d] \ar[l] & E \ar[d]^{\pi} \\
          A & A' \ar@{->>}[l]^p \ar[r] & U }
\end{displaymath}
where $p$ is epi, as indicated.
\end{definition}

Representability formulates the existence of a weak version of a universe. In this paper, classes of small maps will always be assumed to be representable. The map $\pi$ in the definition of representability is often called the universal small map, even though it is not unique (not even up to isomorphism). Representability has the consequence that, in the internal logic of $\bf E$, a map $f: B \to A$ belongs to $S$ iff it holds that
\[ \forall a \in A \exists u \in U: B_a \cong E_u \] 

The axioms for a class of small maps that we have given so far form the basic definition. The definition can be extended by adding various choice or collection principles. As a matter of fact, we will frequently assume that a class of small maps satisfies the collection axiom in the sense of Joyal and Moerdijk in \cite{JoMoe:AST}:
\begin{itemize}
\item[(CA)] For any small map $f: A \to X$ and epi $C \to A$, there exists a quasi-pullback of the form
\begin{displaymath}
\xymatrix{ B \ar[r] \ar[d]_g & C \ar@{->>}[r] & A \ar[d]^f \\
           Y \ar@{->>}[rr] & & X }
\end{displaymath}
where $Y \to X$ is epi and $g: B \to Y$ is small.
\end{itemize}
In \cite{MoerPalm:CST}, Moerdijk and Palmgren work with a much stronger axiom: what they call the axiom of multiple choice AMC (for a precise formulation, see loc.cit.). One of the purposes of this paper is to eliminate the need for this axiom.

As discussed in \cite{MoerPalm:CST}, the collection axiom can be reformulated using the notion of a \emph{collection map}. Informally, a map $g: D \to C$ in $\bf E$ is a collection map, whenever it is true (in the internal logic of $\bf E$), that for any map $\xymatrixcolsep={1.8em}\xymatrix{f: F \ar@{->>}[r] & D_c}$ covering some fibre of $g$, there is another fibre $D_{c'}$ covering $D_c$ via a map $\xymatrixcolsep={1.8em}\xymatrix{p: D_{c'} \ar@{->>}[r] & D_c}$ which factors through $f$. Diagrammatically, one can express this by asking that for any map $c: T \to C$ and any epi $E \to T \times_C D$ there is a diagram of the form
\begin{displaymath}
\xymatrix{ D \ar[d] & D \times_C T' \ar[d] \ar[l] \ar[r] & E \ar@{->>}[r] & T \times_C D \ar[d] \ar[r] & D \ar[d] \\
C & T' \ar[l] \ar@{->>}[rr] & & T \ar[r] & C }
\end{displaymath}
where the middle square is a quasi-pullback with an epi on the bottom, while the two outer squares are pullbacks. A map $g: D \to C$ over $A$ is a \emph{collection map over $A$}, if it is a collection map in ${\bf E}/A$.

The collection axiom is now equivalent to stating that the universal small map $\pi: E \to U$ is a collection map. (This is imprecise, but in a harmless way: if one universal small map is a collection map, they all are.)

We will need the following variation on the notion of a collection map. A span $(g, h)$
\begin{displaymath}
\xymatrix{ D \ar[r]^h \ar[d]_g & B  \\
           C           &  }
\end{displaymath} 
is called a \emph{collection span}, when, in the internal logic, it holds that for any map $\xymatrixcolsep={1.8em}\xymatrix{f: F \ar@{->>}[r] & D_c}$ covering some fibre of $g$, there is another fibre $D_{c'}$ of $g$ covering $D_c$ via a map $\xymatrixcolsep={1.8em}\xymatrix{p: D_{c'} \ar@{->>}[r] & D_c}$ over $B$ which factors through $f$.

We return to the issue of the various possible notions of a predicative topos. In \cite{MoerPalm:CST}, Moerdijk and Palmgren take what they call ``stratified pseudotoposes'' as their notion of predicative topos.
\begin{definition}
A $\Pi W$-pretopos $\bf E$ is called a \emph{stratified pseudotopos}, if it is equip\-ped with a collection (``hierarchy'') of classes of small maps $(S_n)_{n \in \mathbb{N}}$ such that $S_n \subseteq S_{n+1}$ for every $n$, every map in $\bf E$ is contained in some $S_n$, and a specified universal small map $\pi_n: E_n \to U_n$ for $S_n$ has $S_{n+1}$-small codomain.
\end{definition}
Again, it is possible to strengthen the notion of a stratified pseudotopos by assuming that every $S_n$ satisfies some additional choice or collection principles. In \cite{MoerPalm:CST}, it is assumed that every $S_n$ satisfies AMC. Here, we will show how to work with CA instead.

My preferred notion of a predicative topos is even weaker. A $\Pi W$-pretopos $\bf E$ satisfies the universe operator axiom, whenever
\begin{itemize}
\item[(UO)]
Every map $f: B \to A$ is contained in a class of small maps satisfying the collection axiom.
\end{itemize}
This axiom is inspired by Palmgren's notion of a universe operator (see \cite{Palm:InfS} and \cite{Palm:Univ}). His idea was to add to type theory for every dependent type
\begin{quote}
$B(a)$ Set $[a:A]$
\end{quote}
a universe closed under all the type theoretic operations and containing the types $B(a)$ for every $a \in A$ (but not necessarily $A$ itself). Actually, when intensional Martin-L\"of type theory is extended with such a universe operator, the category of setoids will satisfy UO.

This means we have (at least) five notions of a predicative topos. We list them for future reference:
\begin{enumerate} \label{list}
\item A $\Pi W$-pretopos.
\item A stratified pseudotopos without any choice or collection.
\item A stratified pseudotopos in which the classes of small maps satisfy CA.
\item A stratified pseudotopos in which the classes of small maps satisfy AMC.
\item A $\Pi W$-pretopos satisfying UO.\footnote{As for implications between these various notions, the least obvious ones are (4) $\Rightarrow$ (3) $\Rightarrow$ (5). (2) and (5) ought to be incomparable.}
\end{enumerate}

\section{Sites}

\subsection{Different notions of sites}

Here we will give precise definitions of various notions of internal sites. In the following sections we will prove certain equivalences between them in the context of $\Pi W$-pretoposes. 

The basic categorical structure of an (internal) site consists of an internal category $\bf C$ together with a collection of covering families Cov$(C)$ for every object $C$ of $\bf C$. This is formalized by a commutative square of the form
\begin{equation} \label{site}
\begin{gathered}
\xymatrix{ *+\txt{ \underline{Cov} } \ar[r]^m \ar[d]_{\phi} & C_1 \ar[d]^{\txt{cod}} \\
*+\txt{Cov} \ar[r] & C_0 }
\end{gathered}
\end{equation}
(As usual, $C_1$ is the object of arrows and $C_0$ the object of objects of the internal category $\bf C$. And cod is of course the codomain map.)
So any $U \in$ Cov$(C)$ gives rise to an indexing set \underline{Cov}$_U$, indexing a family of arrows all with codomain $C$. Such a covering family $U$ will therefore typically be denoted by $( \alpha_i: C_i \to C \, | \, i \in I )$, where $I$ is the indexing set. 

For a \emph{site}, the following axiom should hold in the internal logic:
\begin{description}
\item[(C)] For any covering family $( \alpha_i: C_i \to C \, | \, i \in I )$ of $C$ and any arrow $f: D \to C$, there exists a covering family $( \beta_j: D_j \to D \, | \, j \in J )$ such that every composite $f\beta_j$ factors through some $\alpha_i$.
\end{description}
A \emph{Grothendieck site} satisfies the following additional requirements:
\begin{description}
\item[(M)] For any object $C$, there is a covering family $U \in \mbox{Cov}(C)$ such that $(1_C: C \to C) \in U$.
\item[(L)] Whenever there is a covering family $( \alpha_i: C_i \to C \, | \, i \in I )$ of $C$ and families $(\beta_{ij}: \xymatrix{ C_{ij} \ar[r] & C_{i}} \, | \, j \in I_i)$ covering $C_i$ for every $i \in I$, there is a family $(\gamma_l: D_l \to C \, | \, l \in L)$ such that for every $\gamma_l$ factors through some $\alpha_i\beta_{ij}$.
\end{description}

A site will be called \emph{strong}, if it satisfies condition (C) in the strong form that the ``pullback'' $(\beta_j \, | \, j \in J)$ is given (externally) as a function of $f$ and $(\alpha_i \, | \, i \in I)$. A Grothendieck site will be called strong, if it satisfies both (C) and (L) in the strong form, i.e., both the ``pullback'' in (C) and the ``composition'' $(\gamma_l: D_l \to C \, | \, l \in L)$ in (L) are given as a function of the initial data. (When the axiom of choice is not externally valid, this is a considerable strengthening of the original definition.) 

In this paper, a pivotal notion is that of a collection site\footnote{The following definition is what was intended, but inaccurately formalized, in \cite{MoerPalm:CST}.}: a site is a \emph{collection site}, when the span $(\phi, m)$ in (\ref{site}) is a collection span over $C_0$.

In the presence of a class of small maps, a site is said to have \emph{small covers}, when all indexing sets are small. Diagrammatically this is expressed by requiring that $\phi$ in (\ref{site}) is small. Finally, a site is called \emph{small}, when every arrow or object in (\ref{site}) is small.

\subsection{Equivalent Grothendieck sites}

As Johnstone explains in \cite{johnstone02a}, in topos theory (M) and (L) are closure conditions that ``might just as well be there'', but are not essential to the notion of a site. This is backed up by the result that in a topos there is for every internal site an equivalent Grothendieck site (equivalent in the sense that it leads to an equivalent category of internal sheaves). The Grothendieck site is inductively generated by closing off the site under the conditions (M) and (L). $\Pi W$-pretoposes contain inductive definitions in the shape of W-types, so one may expect versions of this result to hold in the context of $\Pi W$-pretoposes as well. This section is devoted to the proof that this is in fact the case.
 
\begin{theorem} \label{eqGrsites}
Let $\bf E$ be a $\Pi W$-pretopos and let $\bf C$ be a site in $\bf E$. 
\begin{enumerate}
\item When $\bf C$ is a strong site, there exists an equivalent strong Grothendieck site $\bf D$ in $\bf E$ with the same underlying category. 
\item When $\bf C$ is a collection site, there exists an equivalent collection Grothendieck site $\bf D$ in $\bf E$ with the same underlying category. 
\end{enumerate}
\end{theorem}
\begin{proof}
The construction of $\bf D$ uses the theory of dependent polynomial functors, their initial algebras and its applications, as developed in the work by Gambino and Hyland \cite{GambinoHyland03}. We outline the construction, which is the same both for (1) and (2).

We take the same underlying category $C$. We wish to find a new object of covering families COV over $C_0$ and it should satisfy ($C \in C_0$):
\[ \mbox{COV}_C = 1 + \sum_{U \in \mbox{Cov}(C)} \prod_{i \in \underline{\mbox{Cov}}_U} \mbox{COV}_{\mbox{dom}(m(i))} \]
As shown by Gambino and Hyland, such an object can be constructed by first defining a functor $F: {\bf E}/C_0 \to {\bf E}/C_0$ as follows ($C \in C_0$):
\[ (FX)_C = 1 + \sum_{U \in \mbox{Cov}(C)} \prod_{i \in \underline{\mbox{Cov}}_U} X_{\mbox{dom}(m(i))} \]
This is what they call a dependent polynomial functor, and, in the presence of W-types, these have initial algebras. These algebras are fixed points for the functor and because of their initiality, they allow definition by recursion on their elements. This we use to define the new object \underline{COV} over COV and the new arrow $M$, thereby completing the definition of the site $\bf D$: 
\begin{equation*} 
\begin{gathered}
\xymatrix{ *+\txt{ \underline{COV} } \ar[r]^M \ar[d] & C_1 \ar[d]^{\txt{cod}} \\
*+\txt{COV} \ar[r] & C_0 }
\end{gathered}
\end{equation*}

Elements in COV (over $C \in C_0$) are either * (the unique element of 1) or of the form sup$_U(t)$, where $U \in \mbox{Cov}(C)$ and $t: \underline{\mbox{Cov}}_U \to \mbox{COV}$.
\begin{eqnarray*}
\mbox{\underline{COV}}_* & = & 1 \\
\mbox{\underline{COV}}_{\mbox{sup}_U(t)} & = & \sum_{i \in \underline{\mbox{Cov}}_U} \mbox{\underline{COV}}_{t(i)}
\end{eqnarray*}
The definition of $M$ runs as follows. The unique element of \underline{COV}$_*$ is sent to $1_C: C \to C$. An element $j$ in \underline{COV}$_{\mbox{sup}_U(t)}$ is of the form $(i, k)$, with $i \in \underline{\mbox{Cov}}_U$ and $k \in \mbox{\underline{COV}}_{t(i)}$. This element $(i, k)$ is sent to $m(i) \circ M(k)$. 

We will briefly indicate why the constructed objects have the desired properties. It is easy to see that the covering families are now closed under (M). By induction on the construction of the covering family $(\alpha_i \, | \, i \in I)$, one can see that the covering families are closed under (C) and (L) in their strong form, when they were true in the strong form in the original site, or that they are true (in their normal form), whenever $\bf C$ is a collection site. The proof that any sheaf for $\bf C$ also satisfies the sheaf condition for this covering family relies on a similar proof by induction, as does the proof in case (2) that $\bf D$ inherits the property of being a collection site.

To see that any sheaf for this Grothendieck site is also a sheaf for $\bf C$, take an element $U \in \mbox{Cov}(C)$. Then for every $i \in \underline{\mbox{Cov}}_U$, consider $D_i = \mbox{dom}(\mbox{map}(i))$. The family consisting solely of the identity on $D_i$ is in COV($D_i$). This element in COV($D_i$) is given as a function $t$ of $i \in \underline{\mbox{Cov}}_U$ and we can therefore construct the element $\mbox{sup}_U(t)$ of COV($C$). This covering family consists of the same maps as $U$, therefore a presheaf satisfying the sheaf condition for $\mbox{sup}_U(t)$ also satisfies the sheaf condition for $U$.
\end{proof}
\begin{remark} Although we will not need it, it is good to point out that in case $\bf E$ is equipped with a class of small maps and the original site has small covers or is small, the same will hold for the equivalent Grothendieck site constructed in the proof.
\end{remark}

\subsection{Equivalent collection sites}

Here we want to investigate conditions under which sites in a $\Pi W$-pretopos have equivalent collection sites. Although the following argument is not very difficult, it is a key step in this paper. We fix a $\Pi W$-pretopos $\bf E$.

\begin{lemma} \label{collsp} Suppose $\bf E$ is equipped with a class of small maps satisfying collection. Then any small $f: B \to A$ fits into quasi-pullback diagram 
\begin{displaymath}
\xymatrix{ D \ar@{->>}[r]^h \ar[d]_g & B \ar[d]^f \\
           C \ar@{->>}[r]          & A }
\end{displaymath}
where $(g, h)$ is a collection span over A and $g$ is small.
\end{lemma}
\begin{proof} We use that the universal small map $\pi: E \to U$ is a collection map.

Let $C = \Sigma_{a \in A, u \in U} \{ p:E_u \to B_a \, | \, p \mbox{ is a cover} \}$. The fibre of $g$ above an element $(a, u, p)$ is $E_u$. We leave all the verifications to the reader.
\end{proof}

\begin{proposition} \label{equivcollsite} Suppose $\bf E$ is equipped with a class of small maps satisfying collection. Then for every site $\bf C$ with small covers, there exists an equivalent collection site with small covers and the same underlying category.
\end{proposition} 
\begin{proof}
$\bf C$ is a site with small covers, so in the diagram
\begin{displaymath}
\begin{gathered}
\xymatrix{ *+\txt{ \underline{Cov} } \ar[r]^m \ar[d]_{\phi} & C_1 \ar[d]^{\txt{cod}} \\
*+\txt{Cov} \ar[r] & C_0 }
\end{gathered}
\end{displaymath}
$\phi$ is small. We now apply the previous lemma to $\phi$:
\begin{displaymath}
\begin{gathered} 
\xymatrix{ E \ar@{->>}[r]^h \ar[d]_{\psi} & *+\txt{ \underline{Cov} } \ar[r]^m \ar[d]_{\phi} & C_1 \ar[d]^{\txt{cod}} \\ D \ar@{->>}[r] &
*+\txt{Cov} \ar[r] & C_0 }
\end{gathered}
\end{displaymath}
so the left square is a quasi-pullback, $\psi$ is small and $(\psi, h)$ is a collection span over Cov. Now, the outer rectangle defines an equivalent site with small covers and the same underlying category. It is also collection site, because $(\psi, mh)$ is a collection span over $C_0$.
\end{proof}

\section{Categories of presheaves}

This section is devoted to a proof of the fact that all the notions of a predicative topos contained on the list are closed under taking presheaves. This is new, except for predicative toposes of type 4.

We show this first for the predicative toposes of the simplest kind (i.e., of type 1 on the list).

\begin{theorem} \label{presh1}
If $\bf E$ is a $\Pi W$-pretopos and $\bf C$ is an internal category in $\bf E$, then {\rm PSh}$_{\bf E}({\bf C})$ is a $\Pi W$-pretopos.
\end{theorem}
\begin{proof}
The fact that PSh$_{\bf C}(\bf{E})$ is a locally cartesian closed pretopos with natural number object is well-known. It remains to show that it has W-types. Here we follow \cite{MoerPalm:Wellftr}, pp. 205-8, closely.

With any internal presheaf $P$ in $\bf E$, one can associate the ``underlying set'' $|P|$ given by:
\[ |P| = \{ \, (x, C) \, | \, C \in {\bf C}, x \in P(C) \, \} \]

For a morphism of presheaves $f: B \to A$, and an element $a \in A(C)$, one sets
\[ B_a(D) = \{ \, (\beta: D \to C, b \in B(D) ) \, | \, f(b) = a \cdot \beta \, \} \]
$B_a$ has the structure of a presheaf, when restriction along a morphism $\delta: E \to D$ is defined as:
\[ (\beta, b) \cdot \delta = (\beta\delta, b \cdot \delta) \]

Whenever $X$ is a presheaf, $P_f(X)$ can be written on an object $C$ of $\bf C$ as
\[ P_f(X)(C) = \{(a, t)\,|\, a \in A(C), t: B_a \to X \} \]
where $t$ is a morphism of presheaves. Restriction along a morphism $\alpha: D \to C$ is then given by
\[ (a, t) \cdot \alpha = (a \cdot \alpha, \alpha^{*}(t) ) \]
where $\alpha^{*}(t)(\beta, b) = t(\alpha\beta, b)$.

The presheaf morphism $f$ induces a map 
\[ g: \Sigma_{(a, C) \in |A|} |B_a| \to |A| \]
in $\bf E$ whose fibre over $(a, C)$ is precisely $|B_a|$. The W-type in presheaves will be constructed from the W-type $V$ associated to $g$ in $\bf E$.

This means that every element in $T \in V$ is of the form
\[ T = \sup_{(a, C)} t \]
where $(a, C) \in |A|$ and $t$ is a function $|B_{a}| \to V$. For any such term $T$, one defines its root $\gamma(T)$ to be $C$. If one writes $V(C)$ for the set of terms $T$such that $\gamma(T) = C$, then $V$ has the structure of a presheaf. Restriction along a map $\alpha: C' \to C$ is given by
\[ T \cdot \alpha = \sup_{(a \cdot \alpha, C')} \alpha^{*}(t) \]

In \cite{MoerPalm:Wellftr}, composable and natural terms are defined using an additional axiom allowing transfinite recursion, but our point here is that this axiom is unnecessary. Indeed, composibality and naturality of terms are definable in the internal logic, because the terms with these properties are precisely those for which all the subterms have a certain definable property. We call a term $T$ \emph{composable} if all subterms $\sup_{(a, C)} (t)$ of $T$ have the property that for all $(\beta: D \to C, b) \in B_{a}$,
\[ \gamma(t(\beta, b)) = \mbox{dom}(\beta) \]
A term $T$ is \emph{natural}, if it is composable and if all subterms $\sup_{(a, C)}(t)$ of $T$ have the property that for any $(\beta: D \to C, b) \in B_{a}$ and any $\gamma: E \to D$
\[ t(\beta, b) \cdot \gamma = t(\beta \gamma, b \cdot \gamma) \]
(So $t$ is actually a natural transformation.) These properties are trivially inherited by subterms.

Moerdijk and Palmgren prove that for a natural term $T$ rooted in $C$ and map $\alpha: C' \to C$, the term $T \cdot \alpha$ is also natural (the proof of lemma 5.5 in \cite{MoerPalm:Wellftr} can be copied verbatim). So when $W(C) \subseteq V(C)$ is the collection of natural terms rooted in $C$, $W$ is a subpresheaf of $V$.

We now use the characterization theorem to show that $W$ is a W-type. First of all, there is morphism $S: P_f(W) \to W$, making $W$ into a $P_f$-algebra, because for any $a \in A(C)$ and natural transformation $t: B_a \to W$, one can put
\[ S_C(a, t) = \mbox{sup}_{(a, C)}t \]
$S$ is well-defined and an isomorphism of presheaves, because every natural term $T \in W(C)$ can uniquely be written as $\sup_{(a, C)}t$ for a natural transformation $t: B_a \to W$.

It remains to verify that $W$ has no proper $P_f$-subalgebras. This is easy, because when $K$ is some $P_f$-subalgebra of $W$, then 
\[ L = \{ T \in V \, | \, T \in K(\gamma(T)) \} \]
is a $P_g$-subalgebra of $V$.
\end{proof}

The following proposition takes care of the other types.
\begin{proposition} \label{smallinpresh}
Let $\bf E$ be a $\Pi W$-pretopos with a class of small maps $S$ and $\bf C$ be an internal category in which the codomain map $\mbox{cod}: C_1 \to C_0$ is small. Then {\rm PSh}$_{\bf E}({\bf C})$ inherits a class of small maps, denoted by $T$, by declaring a morphism of presheaves $f: B \to A$ to be $T$-small, whenever
\[ f_C: B(C) \to A(C) \]
is $S$-small for every $C \in C_0$. More formally, if $| \ldots |$ is the forgetful functor {\rm PSh}$_{\bf E}({\bf C}) \to {\bf E}/{C_0}$, a map $f$ is $T$-small, when $\Sigma_{C_0} |f|$ is $S$-small. Moreover, if $S$ satisfies CA or AMC, so does $T$.
\end{proposition}
\begin{proof}
The argument is essentially contained in both \cite{JoMoe:AST} and \cite{MoerPalm:CST} (for the more general case of sheaves), therefore we give only a brief indication of why this result holds.

It is straight forward to see that $T$ is a locally full subcategory, because pullbacks and sums are computed pointwise and the epis in presheaves are precisely those morphisms that are pointwise epic.

Quotients of equivalence relations are also computed pointwise, while exponentials of small objects are constructed using the Yoneda Lemma (see \cite{MacLMoe:Sheaves}, prop. 3.6.1), and are small, because cod is assumed to be small. We leave it to the reader to see that small objects are closed under $\Pi$ and W.

Representability of the class $T$ is proved both in \cite{JoMoe:AST} and \cite{MoerPalm:CST}. Finally, stability of CA can be found in \cite{JoMoe:AST}, while that of AMC can be found in \cite{MoerPalm:CST}.
\end{proof}

\begin{theorem} 
All types of predicative toposes contained in the list are closed under taking presheaves for an internal category.
\end{theorem}
\begin{proof}
For type 1, this statement is precisely theorem \ref{presh1}.

If $\bf E$ is a predicative topos of type 2 and $\bf C$ is an internal category in $\bf E$, let $S_n$ be a class of small maps in the hierarchy such that the codomain map of $\bf C$ is $S_n$-small. Then let $T_m$ be the class of maps in {\rm PSh}$_{\bf E}({\bf C})$ determined by $S_{n+m}$ as in proposition \ref{smallinpresh}. Then {\rm PSh}$_{\bf E}({\bf C})$ is again of type 2. The same argument works for types 3 and 4, because CA and AMC are stable under taking presheaves.

If $\bf E$ is a predicative topos of type 5 and $\bf C$ is an internal category in $\bf E$ and $f$ is an arbitrary map in {\rm PSh}$_{\bf E}({\bf C})$, let $S$ be a class of small maps satisfying collection such that both the codomain map of $\bf C$ and the underlying map $\Sigma_{C_0} |f|$ in $\bf E$ are contained in it. Then the class of small maps $T$ determined by $S$ in {\rm PSh}$_{\bf E}({\bf C})$ satisfies collection and contains $f$.
\end{proof}

\section{Categories of sheaves}

Without too much effort, we can show that the results in the previous section imply the following:

\begin{corollary} \label{immforsh}
If $\bf E$ is a $\Pi W$-pretopos and $\bf C$ is an internal site in $\bf C$, then the categories {\rm Sep}$_{\bf E}({\bf C})$ of separated presheaves and {\rm Sh}$_{\bf E}({\bf C})$ of sheaves are both locally cartesian closed regular categories with W-types.
\end{corollary}
\begin{proof}
Both {\rm Sep}$_{\bf E}({\bf C})$ and {\rm Sh}$_{\bf E}({\bf C})$ have finite limits, because these are computed as in presheaves. {\rm Sep}$_{\bf E}({\bf C})$ is regular, because subobjects of separated presheaves are also separated. To see that {\rm Sh}$_{\bf E}({\bf C})$ is regular, one uses that covers in sheaves are maps that are locally surjective (see \cite{MacLMoe:Sheaves}, p. 143).

It is well-known that exponentials (and the $\Pi$-functors generally) are computed in sheaves as in presheaves. That the same is true for W-types is proposition 5.7 in \cite{MoerPalm:Wellftr}. Similar statements hold for separated presheaves.
\end{proof}

Unfortunately, it appears that one cannot do better: it seems impossible to show that Sh$_{\bf E}({\bf C})$ is a pretopos (more specifically, that it has finite sums and is exact; of course, {\rm Sep}$_{\bf E}({\bf C})$ has finite disjoint sums). This means that we cannot show that predicative toposes of type 1 are closed under taking sheaves for an internal site. In fact, we have similar difficulties for predicative toposes of type 2.

In particular, we are unable to construct a sheafification functor (a left adjoint to the inclusion of sheaves in presheaves) in general. We do, however, have the following result by Moerdijk and Palmgren (lemma 8.1 in \cite{MoerPalm:CST}), which seems to be the best one can say:
\begin{lemma} \label{sheafification}
If a $\Pi W$-pretopos $\bf E$ contains a Grothendieck collection site $\bf C$, the inclusion of sheaves into presheaves $\bf i$ has a left adjoint, the associated sheaf functor
\[ {\bf a}: {\rm PSh}_{\bf E}({\bf C}) \to {\rm Sh}_{\bf E}({\bf C}) \]
\end{lemma}
This means that only when the internal site $\bf C$ is a collection site, will we know that ${\rm Sh}_{\bf E}({\bf C})$ is a $\Pi W$-pretopos. The remainder of the section establishes stability under sheaves of those types of predicative toposes in which every site is equivalent to a collection site. But first we collect the results we have obtained so far in:

\begin{theorem} \label{improve}
If $\bf E$ is a $\Pi W$-pretopos equipped with a class of small maps $S$ satisfying collection and containing an internal site $\bf C$ with small covers, then {\rm Sh}$_{\bf E}({\bf C})$ is a $\Pi W$-pretopos.
\end{theorem}
\begin{proof} By \ref{equivcollsite} and \ref{eqGrsites}, there exists in $\bf E$ an Grothendieck collection site equivalent to $\bf C$ with the same underlying category. Using the left adjoint from \ref{sheafification}, we show that Sh$_{\bf E}({\bf C})$ has finite sums and quotients of equivalence relations (for they can now be calculated in presheaves and then be sheafified). So it is a $\Pi W$-pretopos by \ref{immforsh}.
\end{proof}

\noindent
Observe that this is an improvement over \cite{MoerPalm:CST}, because there this result depends on AMC.

To show that a class of predicative toposes is stable under sheafification, we again need a proposition of the following type:
\begin{proposition} \label{smallinsh}
Let $\bf E$ be a $\Pi W$-pretopos with a class of small maps $S$ and a small Grothendieck collection site $\bf C$. Then {\rm Sh}$_{\bf E}({\bf C})$ inherits a class of small maps, denoted by $T$, by declaring a morphism of sheaves $f: B \to A$ to be $T$-small, whenever
\[ f_C: B(C) \to A(C) \]
is $S$-small for every $C \in C_0$. More formally, when $| \ldots |$ is the forgetful functor {\rm Sh}$_{\bf E}({\bf C}) \to {\bf E}/{C_0}$, a map $f$ is $T$-small, whenever $\Sigma_{C_0} |f|$ is $S$-small. Moreover, if $S$ satisfies CA or AMC, so does $T$.
\end{proposition}
\begin{proof}
The argument is similar to that in \ref{smallinpresh} and again the main points already appeared in \cite{JoMoe:AST} and \cite{MoerPalm:CST}. We therefore refrain from giving a proof.
\end{proof}

\begin{theorem} Predicative toposes of type 5 are closed under taking sheaves for an internal site.
\end{theorem}
\begin{proof}
Let $\bf E$ be a $\Pi W$-pretopos satisfying UO and $\bf C$ be an internal site. We can find a class of small maps $R$ satisfying collection such that $\bf C$ has $R$-small covers, so \ref{improve} implies that {\rm Sh}$_{\bf E}({\bf C})$ is a $\Pi W$-pretopos.

Let $f: B \to A$ be an arbitrary map in {\rm Sh}$_{\bf E}({\bf C})$. Now find a class of small maps $S$ satisfying collection in $\bf E$ such that both the Grothendieck collection site equivalent to $\bf C$ and the map $\Sigma_{C_0} |f|$ in $\bf E$ are $S$-small. Let $T$ be the class of maps in {\rm Sh}$_{\bf E}({\bf C})$ determined by $S$ as in \ref{smallinsh}. Then $T$ is a class of small maps satisfying collection and $f$ is $T$-small.
\end{proof}

\begin{theorem} Predicative toposes of types 3 and 4 are closed under taking sheaves for an internal site.
\end{theorem}
\begin{proof}
Let $\bf E$ be a stratified pseudotopos in which the classes of small maps satisfy either CA or AMC and let $\bf C$ be an internal site in $\bf E$. Then $\bf E$ satisfies UO, so, by theorem \ref{improve}, {\rm Sh}$_{\bf E}({\bf C})$ is a $\Pi W$-pretopos.

The Grothendieck collection site $\bf D$ equivalent to $\bf C$ will be $S_n$-small for some class of small maps $S_n$ in the hierachy. So when $T_m$ is the class of small maps in {\rm Sh}$_{\bf E}({\bf C})$ determined by $S_{m+n}$ as in \ref{smallinsh}, the category of sheaves will be equipped with a hierarchy of small maps making it into a predicative topos of type 3 or 4.
\end{proof}

\bibliographystyle{plain}
\bibliography{sheaves}

\end{document}